# Frequentist statistics as a theory of inductive inference

Deborah G. Mayo[1] and D. R. Cox[2]

*Virginia Tech and Nuffield College, Oxford*

**Abstract:** After some general remarks about the interrelation between philosophical and statistical thinking, the discussion centres largely on significance tests. These are defined as the calculation of *p*-values rather than as formal procedures for "acceptance" and "rejection." A number of types of null hypothesis are described and a principle for evidential interpretation set out governing the implications of *p*-values in the specific circumstances of each application, as contrasted with a long-run interpretation. A variety of more complicated situations are discussed in which modification of the simple *p*-value may be essential.

## 1. Statistics and inductive philosophy

### 1.1. What is the Philosophy of Statistics?

The philosophical foundations of statistics may be regarded as the study of the epistemological, conceptual and logical problems revolving around the use and interpretation of statistical methods, broadly conceived. As with other domains of philosophy of science, work in statistical science progresses largely without worrying about "philosophical foundations". Nevertheless, even in statistical practice, debates about the different approaches to statistical analysis may influence and be influenced by general issues of the nature of inductive-statistical inference, and thus are concerned with foundational or philosophical matters. Even those who are largely concerned with applications are often interested in identifying general principles that underlie and justify the procedures they have come to value on relatively pragmatic grounds. At one level of analysis at least, statisticians and philosophers of science ask many of the same questions.

- What should be observed and what may justifiably be inferred from the resulting data?
- How well do data confirm or fit a model?
- What is a good test?
- Does failure to reject a hypothesis *H* constitute evidence "confirming" *H*?
- How can it be determined whether an apparent anomaly is genuine? How can blame for an anomaly be assigned correctly?
- Is it relevant to the relation between data and a hypothesis if looking at the data influences the hypothesis to be examined?
- How can spurious relationships be distinguished from genuine regularities?

---

[1]Department of Philosophy and Economics, Virginia Tech, Blacksburg, VA 24061-0126, e-mail: mayod@vt.edu

[2]Nuffield College, Oxford OX1 1NF, UK, e-mail: david.cox@nuffield.ox.ac.uk

77

- How can a causal explanation and hypothesis be justified and tested?
- How can the gap between available data and theoretical claims be bridged reliably?

That these very general questions are entwined with long standing debates in philosophy of science helps explain why the field of statistics tends to cross over, either explicitly or implicitly, into philosophical territory. Some may even regard statistics as a kind of "applied philosophy of science" (Fisher [10]; Kempthorne [13]), and statistical theory as a kind of "applied philosophy of inductive inference". As Lehmann [15] has emphasized, Neyman regarded his work not only as a contribution to statistics but also to inductive philosophy. A core question that permeates "inductive philosophy" both in statistics and philosophy is: What is the nature and role of probabilistic concepts, methods, and models in making inferences in the face of limited data, uncertainty and error?

Given the occasion of our contribution, a session on philosophy of statistics for the second Lehmann symposium, we take as our springboard the recommendation of Neyman ([22], p. 17) that we view statistical theory as essentially a "Frequentist Theory of Inductive Inference". The question then arises as to what conception(s) of inductive inference would allow this. Whether or not this is the only or even the most satisfactory account of inductive inference, it is interesting to explore how much progress towards an account of inductive inference, as opposed to inductive behavior, one might get from frequentist statistics (with a focus on testing and associated methods). These methods are, after all, often used for inferential ends, to learn about aspects of the underlying data generating mechanism, and much confusion and criticism (e.g., as to whether and why error rates are to be adjusted) could be avoided if there was greater clarity on the roles in inference of hypothetical error probabilities.

Taking as a backdrop remarks by Fisher [10], Lehmann [15] on Neyman, and by Popper [26] on induction, we consider the roles of significance tests in bridging inductive gaps in traditional hypothetical deductive inference. Our goal is to identify a key principle of evidence by which hypothetical error probabilities may be used for inductive inference from specific data, and to consider how it may direct and justify (a) different uses and interpretations of statistical significance levels in testing a variety of different types of null hypotheses, and (b) when and why "selection effects" need to be taken account of in data dependent statistical testing.

### *1.2. The role of probability in frequentist induction*

The defining feature of an inductive inference is that the premises (evidence statements) can be true while the conclusion inferred may be false without a logical contradiction: the conclusion is "evidence transcending". Probability naturally arises in capturing such evidence transcending inferences, but there is more than one way this can occur. Two distinct philosophical traditions for using probability in inference are summed up by Pearson ([24], p. 228):

"For one school, the degree of confidence in a proposition, a quantity varying with the nature and extent of the evidence, provides the basic notion to which the numerical scale should be adjusted." The other school notes the relevance in ordinary life and in many branches of science of a knowledge of the relative frequency of occurrence of a particular class of events in a series of repetitions, and suggests that "it is through its link with relative frequency that probability has the most direct meaning for the human mind".

Frequentist induction, whatever its form, employs probability in the second manner. For instance, significance testing appeals to probability to characterize the proportion of cases in which a null hypothesis $H_0$ would be rejected in a hypothetical long-run of repeated sampling, an error probability. This difference in the role of probability corresponds to a difference in the form of inference deemed appropriate: The former use of probability traditionally has been tied to the view that a probabilistic account of induction involves quantifying a degree of support or confirmation in claims or hypotheses.

Some followers of the frequentist approach agree, preferring the term "inductive behavior" to describe the role of probability in frequentist statistics. Here the inductive reasoner "decides to infer" the conclusion, and probability quantifies the associated risk of error. The idea that one role of probability arises in science to characterize the "riskiness" or probativeness or severity of the tests to which hypotheses are put is reminiscent of the philosophy of Karl Popper [26]. In particular, Lehmann ([16], p. 32) has noted the temporal and conceptual similarity of the ideas of Popper and Neyman on "finessing" the issue of induction by replacing inductive reasoning with a process of hypothesis testing.

It is true that Popper and Neyman have broadly analogous approaches based on the idea that we can speak of a hypothesis having been well-tested in some sense, quite distinct from its being accorded a degree of probability, belief or confirmation; this is "finessing induction". Both also broadly shared the view that in order for data to "confirm" or "corroborate" a hypothesis $H$, that hypothesis would have to have been subjected to a test with high probability or power to have rejected it if false. But despite the close connection of the ideas, there appears to be no reference to Popper in the writings of Neyman (Lehmann [16], p. 3) and the references by Popper to Neyman are scant and scarcely relevant. Moreover, because Popper denied that any inductive claims were justifiable, his philosophy forced him to deny that even the method he espoused (conjecture and refutations) was reliable. Although $H$ might be true, Popper made it clear that he regarded corroboration at most as a report of the past performance of $H$: it warranted no claims about its reliability in future applications. By contrast, a central feature of frequentist statistics is to be able to assess and control the probability that a test would have rejected a hypothesis, if false. These probabilities come from formulating the data generating process in terms of a statistical model.

Neyman throughout his work emphasizes the importance of a probabilistic model of the system under study and describes frequentist statistics as modelling the phenomenon of the stability of relative frequencies of results of repeated "trials", granting that there are other possibilities concerned with modelling psychological phenomena connected with intensities of belief, or with readiness to bet specified sums, etc. citing Carnap [2], de Finetti [8] and Savage [27]. In particular Neyman criticized the view of "frequentist" inference taken by Carnap for overlooking the key role of the stochastic model of the phenomenon studied. Statistical work related to the inductive philosophy of Carnap [2] is that of Keynes [14] and, with a more immediate impact on statistical applications, Jeffreys [12].

### *1.3. Induction and hypothetical-deductive inference*

While "hypothetical-deductive inference" may be thought to "finesse" induction, in fact inductive inferences occur throughout empirical testing. Statistical testing ideas may be seen to fill these inductive gaps: If the hypothesis were deterministic we could find a relevant function of the data whose value (i) represents the relevant feature under test and (ii) can be predicted by the hypothesis. We calculate the function and then see whether the data agree or disagree with the prediction. If the data conflict with the prediction, then either the hypothesis is in error or some auxiliary or other background factor may be blamed for the anomaly (Duhem's problem).

Statistical considerations enter in two ways. If $H$ is a statistical hypothesis, then usually no outcome strictly contradicts it. There are major problems involved in regarding data as inconsistent with $H$ merely because they are highly improbable; all individual outcomes described in detail may have very small probabilities. Rather the issue, essentially following Popper ([26], pp. 86, 203), is whether the possibly anomalous outcome represents some systematic and reproducible effect.

The focus on falsification by Popper as the goal of tests, and falsification as the defining criterion for a scientific theory or hypothesis, clearly is strongly redolent of Fisher's thinking. While evidence of direct influence is virtually absent, the views of Popper agree with the statement by Fisher ([9], p. 16) that every experiment may be said to exist only in order to give the facts the chance of disproving the null hypothesis. However, because Popper's position denies ever having grounds for inference about reliability, he denies that we can ever have grounds for inferring reproducible deviations.

The advantage in the modern statistical framework is that the probabilities arise from defining a probability model to represent the phenomenon of interest. Had Popper made use of the statistical testing ideas being developed at around the same time, he might have been able to substantiate his account of falsification.

The second issue concerns the problem of how to reason when the data "agree" with the prediction. The argument from $H$ entails data $y$, and that $y$ is observed, to the inference that $H$ is correct is, of course, deductively invalid. A central problem for an inductive account is to be able nevertheless to warrant inferring $H$ in some sense. However, the classical problem, even in deterministic cases, is that many rival hypotheses (some would say infinitely many) would also predict $y$, and thus would pass as well as $H$. In order for a test to be probative, one wants the prediction from $H$ to be something that at the same time is in some sense very surprising and not easily accounted for were $H$ false and important rivals to $H$ correct. We now consider how the gaps in inductive testing may bridged by a specific kind of statistical procedure, the significance test.

## 2. Statistical significance tests

Although the statistical significance test has been encircled by controversies for over 50 years, and has been mired in misunderstandings in the literature, it illustrates in simple form a number of key features of the perspective on frequentist induction that we are considering. See for example Morrison and Henkel [21] and Gibbons and Pratt [11]. So far as possible, we begin with the core elements of significance testing in a version very strongly related to but in some respects different from both Fisherian and Neyman-Pearson approaches, at least as usually formulated.

### 2.1. General remarks and definition

We suppose that we have empirical data denoted collectively by $y$ and that we treat these as observed values of a random variable $Y$. We regard $y$ as of interest only in so far as it provides information about the probability distribution of $Y$ as defined by the relevant statistical model. This probability distribution is to be regarded as an often somewhat abstract and certainly idealized representation of the underlying data generating process. Next we have a hypothesis about the probability distribution, sometimes called the hypothesis under test but more often conventionally called the null hypothesis and denoted by $H_0$. We shall later set out a number of quite different types of null hypotheses but for the moment we distinguish between those, sometimes called simple, that completely specify (in principle numerically) the distribution of $Y$ and those, sometimes called composite, that completely specify certain aspects and which leave unspecified other aspects.

In many ways the most elementary, if somewhat hackneyed, example is that $Y$ consists of $n$ independent and identically distributed components normally distributed with unknown mean $\mu$ and possibly unknown standard deviation $\sigma$. A simple hypothesis is obtained if the value of $\sigma$ is known, equal to $\sigma_0$, say, and the null hypothesis is that $\mu = \mu_0$, a given constant. A composite hypothesis in the same context might have $\sigma$ unknown and again specify the value of $\mu$.

Note that in this formulation it is required that some unknown aspect of the distribution, typically one or more unknown parameters, is precisely specified. The hypothesis that, for example, $\mu \leq \mu_0$ is not an acceptable formulation for a null hypothesis in a Fisherian test; while this more general form of null hypothesis is allowed in Neyman-Pearson formulations.

The immediate objective is to test the conformity of the particular data under analysis with $H_0$ in some respect to be specified. To do this we find a function $t = t(y)$ of the data, to be called the test statistic, such that

- the larger the value of $t$ the more inconsistent are the data with $H_0$;
- the corresponding random variable $T = t(Y)$ has a (numerically) known probability distribution when $H_0$ is true.

These two requirements parallel the corresponding deterministic ones. To assess whether there is a genuine discrepancy (or reproducible deviation) from $H_0$ we define the so-called $p$-value corresponding to any $t$ as

$$p = p(t) = P(T \geq t; H_0),$$

regarded as a measure of concordance with $H_0$ in the respect tested. In at least the initial formulation alternative hypotheses lurk in the undergrowth but are not explicitly formulated probabilistically; also there is no question of setting in advance a preassigned threshold value and "rejecting" $H_0$ if and only if $p \leq \alpha$. Moreover, the justification for tests will not be limited to appeals to long run-behavior but will instead identify an inferential or evidential rationale. We now elaborate.

## 2.2. Inductive behavior vs. inductive inference

The reasoning may be regarded as a statistical version of the valid form of argument called in deductive logic *modus tollens*. This infers the denial of a hypothesis $H$ from the combination that $H$ entails $E$, together with the information that $E$ is false. Because there was a high probability $(1-p)$ that a less significant result would have occurred were $H_0$ true, we may justify taking low $p$-values, properly computed, as evidence against $H_0$. Why? There are two main reasons:

Firstly such a rule provides low error rates (i.e., erroneous rejections) in the long run when $H_0$ is true, a behavioristic argument. In line with an error- assessment view of statistics we may give any particular value $p$, say, the following hypothetical interpretation: suppose that we were to treat the data as just decisive evidence against $H_0$. Then in hypothetical repetitions $H_0$ would be rejected in a long-run proportion $p$ of the cases in which it is actually true. However, knowledge of these hypothetical error probabilities may be taken to underwrite a distinct justification.

This is that such a rule provides a way to determine whether a specific data set is evidence of a discrepancy from $H_0$.

In particular, a low $p$-value, so long as it is properly computed, provides evidence of a discrepancy from $H_0$ in the respect examined, while a $p$-value that is not small affords evidence of accordance or consistency with $H_0$ (where this is to be distinguished from positive evidence for $H_0$, as discussed below in Section 2.3). Interest in applications is typically in whether $p$ is in some such range as $p \geq 0.1$ which can be regarded as reasonable accordance with $H_0$ in the respect tested, or whether $p$ is near to such conventional numbers as 0.05, 0.01, 0.001. Typical practice in much applied work is to give the observed value of $p$ in rather approximate form. A small value of $p$ indicates that (i) $H_0$ is false (there is a discrepancy from $H_0$) or (ii) the basis of the statistical test is flawed, often that real errors have been underestimated, for example because of invalid independence assumptions, or (iii) the play of chance has been extreme.

It is part of the object of good study design and choice of method of analysis to avoid (ii) by ensuring that error assessments are relevant.

There is no suggestion whatever that the significance test would typically be the only analysis reported. In fact, a fundamental tenet of the conception of inductive learning most at home with the frequentist philosophy is that inductive inference requires building up incisive arguments and inferences by putting together several different piece-meal results. Although the complexity of the story makes it more difficult to set out neatly, as, for example, if a single algorithm is thought to capture the whole of inductive inference, the payoff is an account that approaches the kind of full-bodied arguments that scientists build up in order to obtain reliable knowledge and understanding of a field.

Amidst the complexity, significance test reasoning reflects a fairly straightforward conception of evaluating evidence anomalous for $H_0$ in a statistical context, the one Popper perhaps had in mind but lacked the tools to implement. The basic idea is that error probabilities may be used to evaluate the "riskiness" of the predictions $H_0$ is required to satisfy, by assessing the reliability with which the test discriminates whether (or not) the actual process giving rise to the data accords with that described in $H_0$. Knowledge of this probative capacity allows determining if there is strong evidence of discrepancy The reasoning is based on the following frequentist principle for identifying whether or not there is evidence against $H_0$:

**FEV (i)** $y$ is (strong) evidence against $H_0$, i.e. (strong) evidence of discrepancy from $H_0$, if and only if, where $H_0$ a correct description of the mechanism generating $y$, then, with high probability, this would have resulted in a less discordant result than is exemplified by $y$.

A corollary of **FEV** is that $y$ is not (strong) evidence against $H_0$, if the probability of a more discordant result is not very low, even if $H_0$ is correct. That is, if there is a moderately high probability of a more discordant result, even were $H_0$ correct, then $H_0$ accords with $y$ in the respect tested.

Somewhat more controversial is the interpretation of a failure to find a small $p$-value; but an adequate construal may be built on the above form of **FEV**.

## *2.3. Failure and confirmation*

The difficulty with regarding a modest value of $p$ as evidence in favour of $H_0$ is that accordance between $H_0$ and $y$ may occur even if rivals to $H_0$ seriously different from $H_0$ are true. This issue is particularly acute when the amount of data is limited. However, sometimes we can find evidence for $H_0$, understood as an assertion that a particular discrepancy, flaw, or error is absent, and we can do this by means of tests that, with high probability, would have reported a discrepancy had one been present. As much as Neyman is associated with automatic decision-like techniques, in practice at least, both he and E. S. Pearson regarded the appropriate choice of error probabilities as reflecting the specific context of interest (Neyman[23], Pearson [24]).

There are two different issues involved. One is whether a particular value of $p$ is to be used as a threshold in each application. This is the procedure set out in most if not all formal accounts of Neyman-Pearson theory. The second issue is whether control of long-run error rates is a justification for frequentist tests or whether the ultimate justification of tests lies in their role in interpreting evidence in particular cases. In the account given here, the achieved value of $p$ is reported, at least approximately, and the "accept- reject" account is purely hypothetical to give $p$ an operational interpretation. E. S. Pearson [24] is known to have disassociated himself from a narrow behaviourist interpretation (Mayo [17]). Neyman, at least in his discussion with Carnap (Neyman [23]) seems also to hint at a distinction between behavioural and inferential interpretations.

In an attempt to clarify the nature of frequentist statistics, Neyman in this discussion was concerned with the term "degree of confirmation" used by Carnap. In the context of an example where an optimum test had failed to "reject" $H_0$, Neyman considered whether this "confirmed" $H_0$. He noted that this depends on the meaning of words such as "confirmation" and "confidence" and that in the context where $H_0$ had not been "rejected" it would be "dangerous" to regard this as confirmation of $H_0$ if the test in fact had little chance of detecting an important discrepancy from $H_0$ even if such a discrepancy were present. On the other hand if the test had appreciable power to detect the discrepancy the situation would be "radically different".

Neyman is highlighting an inductive fallacy associated with "negative results", namely that if data $y$ yield a test result that is not statistically significantly different from $H_0$ (e.g., the null hypothesis of 'no effect'), and yet the test has small probability of rejecting $H_0$, even when a serious discrepancy exists, then $y$ is not good evidence for inferring that $H_0$ is confirmed by $y$. One may be confident in the absence of a discrepancy, according to this argument, only if the chance that the test would have correctly detected a discrepancy is high.

Neyman compares this situation with interpretations appropriate for inductive behaviour. Here confirmation and confidence may be used to describe the choice of action, for example refraining from announcing a discovery or the decision to treat $H_0$ as satisfactory. The rationale is the pragmatic behavioristic one of controlling errors in the long-run. This distinction implies that even for Neyman evidence for deciding may require a distinct criterion than evidence for believing; but unfortunately Neyman did not set out the latter explicitly. We propose that the needed evidential principle is an adaption of **FEV**(**i**) for the case of a $p$-value that is not small:

**FEV**(**ii**): A moderate $p$ value is evidence of the absence of a discrepancy $\delta$ from $H_0$, only if there is a high probability the test would have given a worse fit with $H_0$ (i.e., smaller $p$ value) were a discrepancy $\delta$ to exist. **FEV(ii)** especially arises in the context of "embedded" hypotheses (below).

What makes the kind of hypothetical reasoning relevant to the case at hand is not solely or primarily the long-run low error rates associated with using the tool (or test) in this manner; it is rather what those error rates reveal about the data generating source or phenomenon. The error-based calculations provide reassurance that incorrect interpretations of the evidence are being avoided in the particular case. To distinguish between this "evidential" justification of the reasoning of significance tests, and the "behavioristic" one, it may help to consider a very informal example of applying this reasoning "to the specific case". Thus suppose that weight gain is measured by well-calibrated and stable methods, possibly using several measuring instruments and observers and the results show negligible change over a test period of interest. This may be regarded as grounds for inferring that the individual's weight gain is negligible within limits set by the sensitivity of the scales. Why?

While it is true that by following such a procedure in the long run one would rarely report weight gains erroneously, that is not the rationale for the particular inference. The justification is rather that the error probabilistic properties of the weighing procedure reflect what is actually the case in the specific instance. (This should be distinguished from the evidential interpretation of Neyman–Pearson theory suggested by Birnbaum [1], which is not data-dependent.)

The significance test is a measuring device for accordance with a specified hypothesis calibrated, as with measuring devices in general, by its performance in repeated applications, in this case assessed typically theoretically or by simulation. Just as with the use of measuring instruments, applied to a specific case, we employ the performance features to make inferences about aspects of the particular thing that is measured, aspects that the measuring tool is appropriately capable of revealing.

Of course for this to hold the probabilistic long-run calculations must be as relevant as feasible to the case in hand. The implementation of this surfaces in statistical theory in discussions of conditional inference, the choice of appropriate distribution for the evaluation of $p$. Difficulties surrounding this seem more technical than conceptual and will not be dealt with here, except to note that the exercise of applying (or attempting to apply) **FEV** may help to guide the appropriate test specification.

## 3. Types of null hypothesis and their corresponding inductive inferences

In the statistical analysis of scientific and technological data, there is virtually always external information that should enter in reaching conclusions about what the data indicate with respect to the primary question of interest. Typically, these background considerations enter not by a probability assignment but by identifying the question to be asked, designing the study, interpreting the statistical results and relating those inferences to primary scientific ones and using them to extend and support underlying theory. Judgments about what is relevant and informative must be supplied for the tools to be used non-fallaciously and as intended. Nevertheless, there are a cluster of systematic uses that may be set out corresponding to types of test and types of null hypothesis.

*3.1. Types of null hypothesis*

We now describe a number of types of null hypothesis. The discussion amplifies that given by Cox ([4], [5]) and by Cox and Hinkley [6]. Our goal here is not to give a guide for the panoply of contexts a researcher might face, but rather to elucidate some of the different interpretations of test results and the associated *p*-values. In Section 4.3, we consider the deeper interpretation of the corresponding inductive inferences that, in our view, are (and are not) licensed by *p*-value reasoning.

1. *Embedded null hypotheses.* In these problems there is formulated, not only a probability model for the null hypothesis, but also models that represent other possibilities in which the null hypothesis is false and, usually, therefore represent possibilities we would wish to detect if present. Among the number of possible situations, in the most common there is a parametric family of distributions indexed by an unknown parameter $\theta$ partitioned into components $\theta = (\phi, \lambda)$, such that the null hypothesis is that $\phi = \phi_0$, with $\lambda$ an unknown nuisance parameter and, at least in the initial discussion with $\phi$ one-dimensional. Interest focuses on alternatives $\phi > \phi_0$.

This formulation has the technical advantage that it largely determines the appropriate test statistic $t(y)$ by the requirement of producing the most sensitive test possible with the data at hand.

There are two somewhat different versions of the above formulation. In one the full family is a tentative formulation intended not to so much as a possible base for ultimate interpretation but as a device for determining a suitable test statistic. An example is the use of a quadratic model to test adequacy of a linear relation; on the whole polynomial regressions are a poor base for final analysis but very convenient and interpretable for detecting small departures from a given form. In the second case the family is a solid base for interpretation. Confidence intervals for $\phi$ have a reasonable interpretation.

One other possibility, that arises very rarely, is that there is a simple null hypothesis and a single simple alternative, i.e. only two possible distributions are under consideration. If the two hypotheses are considered on an equal basis the analysis is typically better considered as one of hypothetical or actual discrimination, i.e. of determining which one of two (or more, generally a very limited number) of possibilities is appropriate, treating the possibilities on a conceptually equal basis.

There are two broad approaches in this case. One is to use the likelihood ratio as an index of relative fit, possibly in conjunction with an application of Bayes theorem. The other, more in accord with the error probability approach, is to take each model in turn as a null hypothesis and the other as alternative leading to an assessment as to whether the data are in accord with both, one or neither hypothesis. Essentially the same interpretation results by applying **FEV** to this case, when it is framed within a Neyman–Pearson framework.

We can call these three cases those of a formal family of alternatives, of a well-founded family of alternatives and of a family of discrete possibilities.

2. *Dividing null hypotheses.* Quite often, especially but not only in technological applications, the focus of interest concerns a comparison of two or more conditions, processes or treatments with no particular reason for expecting the outcome to be exactly or nearly identical, e.g., compared with a standard a new drug may increase or may decrease survival rates.

One, in effect, combines two tests, the first to examine the possibility that $\mu > \mu_0$,

say, the other for $\mu < \mu_0$. In this case, the two- sided test combines both one-sided tests, each with its own significance level. The significance level is twice the smaller $p$, because of a "selection effect" (Cox and Hinkley [6], p. 106). We return to this issue in Section 4. The null hypothesis of zero difference then divides the possible situations into two qualitatively different regions with respect to the feature tested, those in which one of the treatments is superior to the other and a second in which it is inferior.

3. *Null hypotheses of absence of structure.* In quite a number of relatively empirically conceived investigations in fields without a very firm theory base, data are collected in the hope of finding structure, often in the form of dependencies between features beyond those already known. In epidemiology this takes the form of tests of potential risk factors for a disease of unknown aetiology.

4. *Null hypotheses of model adequacy.* Even in the fully embedded case where there is a full family of distributions under consideration, rich enough potentially to explain the data whether the null hypothesis is true or false, there is the possibility that there are important discrepancies with the model sufficient to justify extension, modification or total replacement of the model used for interpretation. In many fields the initial models used for interpretation are quite tentative; in others, notably in some areas of physics, the models have a quite solid base in theory and extensive experimentation. But in all cases the possibility of model misspecification has to be faced even if only informally.

There is then an uneasy choice between a relatively focused test statistic designed to be sensitive against special kinds of model inadequacy (powerful against specific directions of departure), and so-called omnibus tests that make no strong choices about the nature of departures. Clearly the latter will tend to be insensitive, and often extremely insensitive, against specific alternatives. The two types broadly correspond to chi-squared tests with small and large numbers of degrees of freedom. For the focused test we may either choose a suitable test statistic or, almost equivalently, a notional family of alternatives. For example to examine agreement of $n$ independent observations with a Poisson distribution we might in effect test the agreement of the sample variance with the sample mean by a chi-squared dispersion test (or its exact equivalent) or embed the Poisson distribution in, for example, a negative binomial family.

5. *Substantively-based null hypotheses.* In certain special contexts, null results may indicate substantive evidence for scientific claims in contexts that merit a fifth category. Here, a theory $\mathcal{T}$ for which there is appreciable theoretical and/or empirical evidence predicts that $H_0$ is, at least to a very close approximation, the true situation.

(a) In one version, there may be results apparently anomalous for $\mathcal{T}$, and a test is designed to have ample opportunity to reveal a discordancy with $H_0$ if the anomalous results are genuine.

(b) In a second version a rival theory $\mathcal{T}^*$ predicts a specified discrepancy from $H_0$. and the significance test is designed to discriminate between $\mathcal{T}$ and the rival theory $\mathcal{T}^*$ (in a thus far not tested domain).

For an example of (a) physical theory suggests that because the quantum of energy in nonionizing electro-magnetic fields, such as those from high voltage transmission lines, is much less than is required to break a molecular bond, there should be no carcinogenic effect from exposure to such fields. Thus in a randomized ex-periment in which two groups of mice are under identical conditions except that one group is exposed to such a field, the null hypothesis that the cancer incidence rates in the two groups are identical may well be exactly true and would be a prime focus of interest in analysing the data. Of course the null hypothesis of this general kind does not have to be a model of zero effect; it might refer to agreement with previous well-established empirical findings or theory.

### 3.2. Some general points

We have in the above described essentially one-sided tests. The extension to two-sided tests does involve some issues of definition but we shall not discuss these here.

Several of the types of null hypothesis involve an incomplete probability specification. That is, we may have only the null hypothesis clearly specified. It might be argued that a full probability formulation should always be attempted covering both null and feasible alternative possibilities. This may seem sensible in principle but as a strategy for direct use it is often not feasible; in any case models that would cover all reasonable possibilities would still be incomplete and would tend to make even simple problems complicated with substantial harmful side-effects.

Note, however, that in all the formulations used here some notion of explanations of the data alternative to the null hypothesis is involved by the choice of test statistic; the issue is when this choice is made via an explicit probabilistic formulation. The general principle of evidence **FEV** helps us to see that in specified contexts, the former suffices for carrying out an evidential appraisal (see Section 3.3).

It is, however, sometimes argued that the choice of test statistic can be based on the distribution of the data under the null hypothesis alone, in effect choosing minus the log probability as test statistic, thus summing probabilities over all sample points as or less probable than that observed. While this often leads to sensible results we shall not follow that route here.

### 3.3. Inductive inferences based on outcomes of tests

How does significance test reasoning underwrite inductive inferences or evidential evaluations in the various cases? The hypothetical operational interpretation of the *p*-value is clear but what are the deeper implications either of a modest or of a small value of $p$? These depends strongly both on (i) the type of null hypothesis, and (ii) the nature of the departure or alternative being probed, as well as (iii) whether we are concerned with the interpretation of particular sets of data, as in most detailed statistical work, or whether we are considering a broad model for analysis and interpretation in a field of study. The latter is close to the traditional Neyman-Pearson formulation of fixing a critical level and accepting, in some sense, $H_0$ if $p > \alpha$ and rejecting $H_0$ otherwise. We consider some of the familiar shortcomings of a routine or mechanical use of *p*-values.

### 3.4. The routine-behavior use of p-values

Imagine one sets $\alpha = 0.05$ and that results lead to a publishable paper if and only for the relevant $p$, the data yield $p < 0.05$. The rationale is the behavioristic one outlined earlier. Now the great majority of statistical discussion, going back to Yates [32] and earlier, deplores such an approach, both out of a concern that it encourages mechanical, automatic and unthinking procedures, as well as a desire to emphasize estimation of relevant effects over testing of hypotheses. Indeed a few journals in some fields have in effect banned the use of $p$-values. In others, such as a number of areas of epidemiology, it is conventional to emphasize 95% confidence intervals, as indeed is in line with much mainstream statistical discussion. Of course, this does not free one from needing to give a proper frequentist account of the use and interpretation of confidence levels, which we do not do here (though see Section 3.6).

Nevertheless the relatively mechanical use of $p$-values, while open to parody, is not far from practice in some fields; it does serve as a screening device, recognizing the possibility of error, and decreasing the possibility of the publication of misleading results. A somewhat similar role of tests arises in the work of regulatory agents, in particular the FDA. While requiring studies to show $p$ less than some preassigned level by a preordained test may be inflexible, and the choice of critical level arbitrary, nevertheless such procedures have virtues of impartiality and relative independence from unreasonable manipulation. While adhering to a fixed $p$-value may have the disadvantage of biasing the literature towards positive conclusions, it offers an appealing assurance of some known and desirable long-run properties. They will be seen to be particularly appropriate for Example 3 of Section 4.2.

### *3.5. The inductive-evidence use of p-values*

We now turn to the use of significance tests which, while more common, is at the same time more controversial; namely as one tool to aid the analysis of specific sets of data, and/or base inductive inferences on data. The discussion presupposes that the probability distribution used to assess the $p$-value is as appropriate as possible to the specific data under analysis.

The general frequentist principle for inductive reasoning, **FEV**, or something like it, provides a guide for the appropriate statement about evidence or inference regarding each type of null hypothesis. Much as one makes inferences about changes in body mass based on performance characteristics of various scales, one may make inferences from significance test results by using error rate properties of tests. They indicate the capacity of the particular test to have revealed inconsistencies and discrepancies in the respects probed, and this in turn allows relating $p$-values to hypotheses about the process as statistically modelled. It follows that an adequate frequentist account of inference should strive to supply the information to implement **FEV**.

*Embedded Nulls.* In the case of embedded null hypotheses, it is straightforward to use small $p$-values as evidence of discrepancy from the null in the direction of the alternative. Suppose, however, that the data are found to accord with the null hypothesis ($p$ not small). One may, if it is of interest, regard this as evidence that any discrepancy from the null is less than $\delta$, using the same logic in significance testing. In such cases concordance with the null may provide evidence of the absence of a discrepancy from the null of various sizes, as stipulated in **FEV**(ii).

To infer the absence of a discrepancy from $H_0$ as large as $\delta$ we may examine the probability $\beta(\delta)$ of observing a worse fit with $H_0$ if $\mu = \mu_0 + \delta$. If that probability is near one then, following **FEV**(ii), the data are good evidence that $\mu < \mu_0 + \delta$. Thus $\beta(\delta)$ may be regarded as the stringency or severity with which the test has probed the discrepancy $\delta$; equivalently one might say that $\mu < \mu_0 + \delta$ has passed a severe test (Mayo [17]).

This avoids unwarranted interpretations of consistency with $H_0$ with insensitive tests. Such an assessment is more relevant to specific data than is the notion of power, which is calculated relative to a predesignated critical value beyond which the test "rejects" the null. That is, power appertains to a prespecified rejection region, not to the specific data under analysis.

Although oversensitivity is usually less likely to be a problem, if a test is so sensitive that a $p$-value as or even smaller than the one observed, is probable even when $\mu < \mu_0 + \delta$, then a small value of $p$ is not evidence of departure from $H_0$ in excess of $\delta$.

If there is an explicit family of alternatives, it will be possible to give a set of confidence intervals for the unknown parameter defining $H_0$ and this would give a more extended basis for conclusions about the defining parameter.

*Dividing and absence of structure nulls.* In the case of dividing nulls, discordancy with the null (using the two-sided value of $p$) indicates direction of departure (e.g., which of two treatments is superior); accordance with $H_0$ indicates that these data do not provide adequate evidence even of the direction of any difference. One often hears criticisms that it is pointless to test a null hypothesis known to be false, but even if we do not expect two means, say, to be equal, the test is informative in order to divide the departures into qualitatively different types. The interpretation is analogous when the null hypothesis is one of absence of structure: a modest value of $p$ indicates that the data are insufficiently sensitive to detect structure. If the data are limited this may be no more than a warning against over-interpretation rather than evidence for thinking that indeed there is no structure present. That is because the test may have had little capacity to have detected any structure present. A small value of $p$, however, indicates evidence of a genuine effect; that to look for a substantive interpretation of such an effect would not be intrinsically error-prone.

Analogous reasoning applies when assessments about the probativeness or sensitivity of tests are informal. If the data are so extensive that accordance with the null hypothesis implies the absence of an effect of practical importance, and a reasonably high $p$-value is achieved, then it may be taken as evidence of the absence of an effect of practical importance. Likewise, if the data are of such a limited extent that it can be assumed that data in accord with the null hypothesis are consistent also with departures of scientific importance, then a high $p$-value does not warrant inferring the absence of scientifically important departures from the null hypothesis.

*Nulls of model adequacy.* When null hypotheses are assertions of model adequacy, the interpretation of test results will depend on whether one has a relatively focused test statistic designed to be sensitive against special kinds of model inadequacy, or so called omnibus tests. Concordance with the null in the former case gives evidence of absence of the type of departure that the test is sensitive in detecting, whereas, with the omnibus test, it is less informative. In both types of tests, a small $p$-value is evidence of some departure, but so long as various alternative models could account for the observed violation (i.e., so long as this test had little ability to discriminate between them), these data by themselves may only provide provisional suggestions of alternative models to try.

*Substantive nulls.* In the preceding cases, accordance with a null could at most provide evidence to rule out discrepancies of specified amounts or types, according to the ability of the test to have revealed the discrepancy. More can be said in the case of substantive nulls. If the null hypothesis represents a prediction from some theory being contemplated for general applicability, consistency with the null hypothesis may be regarded as some additional evidence for the theory, especially if the test and data are sufficiently sensitive to exclude major departures from the theory. An aspect is encapsulated in Fisher's aphorism (Cochran [3]) that to help make observational studies more nearly bear a causal interpretation, one should make one's theories elaborate, by which he meant one should plan a variety of tests of different consequences of a theory, to obtain a comprehensive check of its implications. The limited result that one set of data accords with the theory adds one piece to the evidence whose weight stems from accumulating an ability to refute alternative explanations.

In the first type of example under this rubric, there may be apparently anomalous results for a theory or hypothesis $\mathcal{T}$, where $\mathcal{T}$ has successfully passed appreciable theoretical and/or empirical scrutiny. Were the apparently anomalous results for $\mathcal{T}$ genuine, it is expected that $H_0$ will be rejected, so that when it is not, the results are positive evidence against the reality of the anomaly. In a second type of case, one again has a well-tested theory $\mathcal{T}$, and a rival theory $\mathcal{T}^*$ is determined to conflict with $\mathcal{T}$ in a thus far untested domain, with respect to an effect. By identifying the null with the prediction from $\mathcal{T}$, any discrepancies in the direction of $\mathcal{T}^*$ are given a very good chance to be detected, such that, if no significant departure is found, this constitutes evidence for $\mathcal{T}$ in the respect tested.

Although the general theory of relativity, GTR, was not facing anomalies in the 1960s, rivals to the GTR predicted a breakdown of the Weak Equivalence Principle for massive self-gravitating bodies, e.g., the earth-moon system: this effect, called the Nordvedt effect would be 0 for GTR (identified with the null hypothesis) and non-0 for rivals. Measurements of the round trip travel times between the earth and moon (between 1969 and 1975) enabled the existence of such an anomaly for GTR to be probed. Finding no evidence against the null hypothesis set upper bounds to the possible violation of the WEP, and because the tests were sufficiently sensitive, these measurements provided good evidence that the Nordvedt effect is absent, and thus evidence for the null hypothesis (Will [31]). Note that such a negative result does not provide evidence for all of GTR (in all its areas of prediction), but it does provide evidence for its correctness with respect to this effect. The logic is this: theory $\mathcal{T}$ predicts $H_0$ is at least a very close approximation to the true situation; rival theory $\mathcal{T}^*$ predicts a specified discrepancy from $H_0$, and the test has high probability of detecting such a discrepancy from $\mathcal{T}$ were $\mathcal{T}^*$ correct. Detecting no discrepancy is thus evidence for its absence.

### *3.6. Confidence intervals*

As noted above in many problems the provision of confidence intervals, in principle at a range of probability levels, gives the most productive frequentist analysis. If so, then confidence interval analysis should also fall under our general frequentist principle. It does. In one sided testing of $\mu = \mu_0$ against $\mu > \mu_0$, a small $p$-value corresponds to $\mu_0$ being (just) excluded from the corresponding $(1-2p)$ (two-sided) confidence interval (or $1-p$ for the one-sided interval). Were $\mu = \mu_L$, the lower confidence bound, then a less discordant result would occur with high probability $(1-p)$. Thus **FEV** licenses taking this as evidence of inconsistency with $\mu = \mu_L$ (in the positive direction). Moreover, this reasoning shows the advantage of considering several confidence intervals at a range of levels, rather than just reporting whether or not a given parameter value is within the interval at a fixed confidence level.

Neyman developed the theory of confidence intervals *ab initio* i.e. relying only implicitly rather than explicitly on his earlier work with E.S. Pearson on the theory of tests. It is to some extent a matter of presentation whether one regards interval estimation as so different in principle from testing hypotheses that it is best developed separately to preserve the conceptual distinction. On the other hand there are considerable advantages to regarding a confidence limit, interval or region as the set of parameter values consistent with the data at some specified level, as assessed by testing each possible value in turn by some mutually concordant procedures. In particular this approach deals painlessly with confidence intervals that are null or which consist of all possible parameter values, at some specified significance level. Such null or infinite regions simply record that the data are inconsistent with all possible parameter values, or are consistent with all possible values. It is easy to construct examples where these seem entirely appropriate conclusions.

## 4. Some complications: selection effects

The idealized formulation involved in the initial definition of a significance test in principle starts with a hypothesis and a test statistic, then obtains data, then applies the test and looks at the outcome. The hypothetical procedure involved in the definition of the test then matches reasonably closely what was done; the possible outcomes are the different possible values of the specified test statistic. This permits features of the distribution of the test statistic to be relevant for learning about corresponding features of the mechanism generating the data. There are various reasons why the procedure actually followed may be different and we now consider one broad aspect of that.

It often happens that either the null hypothesis or the test statistic are influenced by preliminary inspection of the data, so that the actual procedure generating the final test result is altered. This in turn may alter the capabilities of the test to detect discrepancies from the null hypotheses reliably, calling for adjustments in its error probabilities.

To the extent that $p$ is viewed as an aspect of the logical or mathematical relation between the data and the probability model such preliminary choices are irrelevant. This will not suffice in order to ensure that the $p$-values serve their intended purpose for frequentist inference, whether in behavioral or evidential contexts. To the extent that one wants the error-based calculations that give the test its meaning to be applicable to the tasks of frequentist statistics, the preliminary analysis and choice may be highly relevant.

The general point involved has been discussed extensively in both philosophical and statistical literatures, in the former under such headings as requiring novelty or avoiding *ad hoc* hypotheses, under the latter, as rules against peeking at the data or shopping for significance, and thus requiring selection effects to be taken into account. The general issue is whether the evidential bearing of data $y$ on an inference or hypothesis $H_0$ is altered when $H_0$ has been either constructed or selected for testing in such a way as to result in a specific observed relation between $H_0$ and $y$, whether that is agreement or disagreement. Those who favour logical approaches to confirmation say no (e.g., Mill [20], Keynes [14]), whereas those closer to an error statistical conception say yes (Whewell [30], Pierce [25]). Following the latter philosophy, Popper required that scientists set out in advance what outcomes they would regard as falsifying $H_0$, a requirement that even he came to reject; the entire issue in philosophy remains unresolved (Mayo [17]).

Error statistical considerations allow going further by providing criteria for when various data dependent selections matter and how to take account of their influence on error probabilities. In particular, if the null hypothesis is chosen for testing because the test statistic is large, the probability of finding some such discordance or other may be high even under the null. Thus, following **FEV(i)**, we would not have genuine evidence of discordance with the null, and unless the $p$-value is modified appropriately, the inference would be misleading. To the extent that one wants the error-based calculations that give the test its meaning to supply reassurance that apparent inconsistency in the particular case is genuine and not merely due to chance, adjusting the $p$-value is called for.

Such adjustments often arise in cases involving data dependent selections either in model selection or construction; often the question of adjusting $p$ arises in cases involving multiple hypotheses testing, but it is important not to run cases together simply because there is data dependence or multiple hypothesis testing. We now outline some special cases to bring out the key points in different scenarios. Then we consider whether allowance for selection is called for in each case.

### *4.1. Examples*

**Example 1.** An investigator has, say, 20 independent sets of data, each reporting on different but closely related effects. The investigator does all 20 tests and reports only the smallest $p$, which in fact is about 0.05, and its corresponding null hypothesis. The key points are the independence of the tests and the failure to report the results from insignificant tests.

**Example 2.** A highly idealized version of testing for a DNA match with a given specimen, perhaps of a criminal, is that a search through a data-base of possible matches is done one at a time, checking whether the hypothesis of agreement with the specimen is rejected. Suppose that sensitivity and specificity are both very high. That is, the probabilities of false negatives and false positives are both very small. The first individual, if any, from the data-base for which the hypothesis is rejected is declared to be the true match and the procedure stops there.

**Example 3.** A microarray study examines several thousand genes for potential expression of say a difference between Type 1 and Type 2 disease status. There are thus several thousand hypotheses under investigation in one step, each with its associated null hypothesis.

**Example 4.** To study the dependence of a response or outcome variable $y$ on an explanatory variable $x$ it is intended to use a linear regression analysis of $y$ on $x$. Inspection of the data suggests that it would be better to use the regression of $\log y$ on $\log x$, for example because the relation is more nearly linear or because secondary assumptions, such as constancy of error variance, are more nearly satisfied.

**Example 5.** To study the dependence of a response or outcome variable $y$ on a considerable number of potential explanatory variables $x$, a data-dependent procedure of variable selection is used to obtain a representation which is then fitted by standard methods and relevant hypotheses tested.

**Example 6.** Suppose that preliminary inspection of data suggests some totally unexpected effect or regularity not contemplated at the initial stages. By a formal test the effect is very "highly significant". What is it reasonable to conclude?

*4.2. Need for adjustments for selection*

There is not space to discuss all these examples in depth. A key issue concerns which of these situations need an adjustment for multiple testing or data dependent selection and what that adjustment should be. How does the general conception of the character of a frequentist theory of analysis and interpretation help to guide the answers?

We propose that it does so in the following manner: Firstly it must be considered whether the context is one where the key concern is the control of error rates in a series of applications (behavioristic goal), or whether it is a context of making a specific inductive inference or evaluating specific evidence (inferential goal). The relevant error probabilities may be altered for the former context and not for the latter. Secondly, the relevant sequence of repetitions on which to base frequencies needs to be identified. The general requirement is that we do not report discordance with a null hypothesis by means a procedure that would report discordancies fairly frequently even though the null hypothesis is true. Ascertainment of the relevant hypothetical series on which this error frequency is to be calculated demands consideration of the nature of the problem or inference. More specifically, one must identify the particular obstacles that need to be avoided for a reliable inference in the particular case, and the capacity of the test, as a measuring instrument, to have revealed the presence of the obstacle.

When the goal is appraising specific evidence, our main interest, **FEV** gives some guidance. More specifically the problem arises when data are used to select a hypothesis to test or alter the specification of an underlying model in such a way that **FEV** is either violated or it cannot be determined whether **FEV** is satisfied (Mayo and Kruse [18]).

**Example 1 (Hunting for statistical significance).** The test procedure is very different from the case in which the single null found statistically significant was preset as the hypothesis to test, perhaps it is $H_{0,13}$ ,the 13th null hypothesis out of the 20. In Example 1, the possible results are the possible statistically significant factors that might be found to show a "calculated" statistical significant departure from the null. Hence the type 1 error probability is the probability of finding at least one such significant difference out of 20, even though the global null is true (i.e., all twenty observed differences are due to chance). The probability that this procedure yields an erroneous rejection differs from, and will be much greater than, 0.05 (and is approximately 0.64). There are different, and indeed many more, ways one can err in this example than when one null is prespecified, and this is reflected in the adjusted *p*-value.

This much is well known, but should this influence the interpretation of the result in a context of inductive inference? According to **FEV** it should. However the concern is not the avoidance of often announcing genuine effects erroneously in a series, the concern is that this test performs poorly as a tool for discriminating genuine from chance effects in this particular case. Because at least one such impressive departure, we know, is common even if all are due to chance, the test has scarcely reassured us that it has done a good job of avoiding such a mistake in this case. Even if there are other grounds for believing the genuineness of the one effect that is found, we deny that this test alone has supplied such evidence.

Frequentist calculations serve to examine the particular case, we have been saying, by characterizing the capability of tests to have uncovered mistakes in inference, and on those grounds, the "hunting procedure" has low capacity to have alerted us to, in effect, temper our enthusiasm, even where such tempering is warranted. If, on the other hand, one adjusts the *p*-value to reflect the overall error rate, the test again becomes a tool that serves this purpose.

Example 1 may be contrasted to a standard factorial experiment set up to investigate the effects of several explanatory variables simultaneously. Here there are a number of distinct questions, each with its associated hypothesis and each with its associated *p*-value. That we address the questions via the same set of data rather than via separate sets of data is in a sense a technical accident. Each *p* is correctly interpreted in the context of its own question. Difficulties arise for particular inferences only if we in effect throw away many of the questions and concentrate only on one, or more generally a small number, chosen just because they have the smallest *p*. For then we have altered the capacity of the test to have alerted us, by means of a correctly computed *p*-value, whether we have evidence for the inference of interest.

**Example 2 (Explaining a known effect by eliminative induction).** Example 2 is superficially similar to Example 1, finding a DNA match being somewhat akin to finding a statistically significant departure from a null hypothesis: one searches through data and concentrates on the one case where a "match" with the criminal's DNA is found, ignoring the non-matches. If one adjusts for "hunting" in Example 1, shouldn't one do so in broadly the same way in Example 2? No.

In Example 1 the concern is that of inferring a genuine, "reproducible" effect, when in fact no such effect exists; in Example 2, there is a known effect or specific event, the criminal's DNA, and reliable procedures are used to track down the specific cause or source (as conveyed by the low "erroneous-match" rate.) The probability is high that we would not obtain a match with person $i$, if $i$ were not the criminal; so, by **FEV**, finding the match is, at a qualitative level, good evidence that $i$ is the criminal. Moreover, each non-match found, by the stipulations of the example, virtually excludes that person; thus, the more such negative results the stronger is the evidence when a match is finally found. The more negative results found, the more the inferred "match" is fortified; whereas in Example 1 this is not so.

Because at most one null hypothesis of innocence is false, evidence of innocence on one individual increases, even if only slightly, the chance of guilt of another. An assessment of error rates is certainly possible once the sampling procedure for testing is specified. Details will not be given here.

A broadly analogous situation concerns the anomaly of the orbit of Mercury: the numerous failed attempts to provide a Newtonian interpretation made it all the more impressive when Einstein's theory was found to predict the anomalous results precisely and without any *ad hoc* adjustments.

**Example 3 (Micro-array data).** In the analysis of micro-array data, a reasonable starting assumption is that a very large number of null hypotheses are being tested and that some fairly small proportion of them are (strictly) false, a global null hypothesis of no real effects at all often being implausible. The problem is then one of selecting the sites where an effect can be regarded as established. Here, the need for an adjustment for multiple testing is warranted mainly by a pragmatic concern to avoid "too much noise in the network". The main interest is in how best to adjust error rates to indicate most effectively the gene hypotheses worth following up. An error-based analysis of the issues is then via the false-discovery rate, i.e. essentially the long run proportion of sites selected as positive in which no effect is present. An alternative formulation is via an empirical Bayes model and the conclusions from this can be linked to the false discovery rate. The latter method may be preferable because an error rate specific to each selected gene may be found; the evidence in some cases is likely to be much stronger than in others and this distinction is blurred in an overall false-discovery rate. See Shaffer [28] for a systematic review.

**Example 4 (Redefining the test).** If tests are run with different specifications, and the one giving the more extreme statistical significance is chosen, then adjustment for selection is required, although it may be difficult to ascertain the precise adjustment. By allowing the result to influence the choice of specification, one is altering the procedure giving rise to the $p$-value, and this may be unacceptable. While the substantive issue and hypothesis remain unchanged the precise specification of the probability model has been guided by preliminary analysis of the data in such a way as to alter the stochastic mechanism actually responsible for the test outcome.

An analogy might be testing a sharpshooter's ability by having him shoot and then drawing a bull's-eye around his results so as to yield the highest number of bull's-eyes, the so-called principle of the Texas marksman. The skill that one is allegedly testing and making inferences about is his ability to shoot when the target is given and fixed, while that is not the skill actually responsible for the resulting high score.

By contrast, if the choice of specification is guided not by considerations of the statistical significance of departure from the null hypothesis, but rather because the data indicates the need to allow for changes to achieve linearity or constancy of error variance, no allowance for selection seems needed. Quite the contrary: choosing the more empirically adequate specification gives reassurance that the calculated $p$-value is relevant for interpreting the evidence reliably. (Mayo and Spanos [19]). This might be justified more formally by regarding the specification choice as an informal maximum likelihood analysis, maximizing over a parameter orthogonal to those specifying the null hypothesis of interest.

**Example 5 (Data mining).** This example is analogous to Example 1, although how to make the adjustment for selection may not be clear because the procedure used in variable selection may be tortuous. Here too, the difficulties of selective reporting are bypassed by specifying all those reasonably simple models that are consistent with the data rather than by choosing only one model (Cox and Snell [7]). The difficulties of implementing such a strategy are partly computational rather than conceptual. Examples of this sort are important in much relatively elaborate statistical analysis in that series of very informally specified choices may be made about the model formulation best for analysis and interpretation (Spanos [29]).

**Example 6 (The totally unexpected effect).** This raises major problems. In laboratory sciences with data obtainable reasonably rapidly, an attempt to obtain independent replication of the conclusions would be virtually obligatory. In other contexts a search for other data bearing on the issue would be needed. High statistical significance on its own would be very difficult to interpret, essentially because selection has taken place and it is typically hard or impossible to specify with any realism the set over which selection has occurred. The considerations discussed in Examples 1-5, however, may give guidance. If, for example, the situation is as in Example 2 (explaining a known effect) the source may be reliably identified in a procedure that fortifies, rather than detracts from, the evidence. In a case akin to Example 1, there is a selection effect, but it is reasonably clear what is the set of possibilities over which this selection has taken place, allowing correction of the $p$-value. In other examples, there is a selection effect, but it may not be clear how to make the correction. In short, it would be very unwise to dismiss the possibility of learning from data something new in a totally unanticipated direction, but one must discriminate the contexts in order to gain guidance for what further analysis, if any, might be required.

## 5. Concluding remarks

We have argued that error probabilities in frequentist tests may be used to evaluate the reliability or capacity with which the test discriminates whether or not the actual process giving rise to data is in accordance with that described in $H_0$. Knowledge of this probative capacity allows determination of whether there is strong evidence against $H_0$ based on the frequentist principle we set out **FEV**. What makes the kind of hypothetical reasoning relevant to the case at hand is not the long-run low error rates associated with using the tool (or test) in this manner; it is rather what those error rates reveal about the data generating source or phenomenon. We have not attempted to address the relation between the frequentist and Bayesian analyses of what may appear to be very similar issues. A fundamental tenet of the conception of inductive learning most at home with the frequentist philosophy is that inductive inference requires building up incisive arguments and inferences by putting together several different piece-meal results; we have set out considerations to guide these pieces. Although the complexity of the issues makes it more difficult to set out neatly, as, for example, one could by imagining that a single algorithm encompasses the whole of inductive inference, the payoff is an account that approaches the kind of arguments that scientists build up in order to obtain reliable knowledge and understanding of a field.